\documentclass[10pt,oneside,leqno]{amsart}
\usepackage{amsxtra}
\usepackage{amsopn}
\usepackage{amsmath,amsthm,amssymb}
\usepackage{amscd}
\usepackage{amsfonts}
\usepackage{latexsym}
\usepackage{verbatim}

\theoremstyle{plain}
\newtheorem{theorem}{Theorem}[section]
\newtheorem{definition}[theorem]{Definition}
\newtheorem{lemma}[theorem]{Lemma}
\newtheorem{proposition}[theorem]{Proposition}
\newtheorem{corollary}[theorem]{Corollary}
\newtheorem{remark}[theorem]{Remark}
\newtheorem{example}[theorem]{Example}
\newtheorem{question}[theorem]{QUESTION}
\newtheorem{remark-question}[section]{Remark-Question}
\newtheorem{conjecture}[section]{Conjecture}

\newcommand{\Hol}{{\rm Hol}}

\newcommand\SU{{\rm SU}}

\newcommand\frg{{\mathfrak g}}
\newcommand\frh{{\mathfrak h}}

\newcommand{\la}{\langle}

\newcommand{\ra}{\rangle}

\begin{document}
\title[Balanced Hermitian metrics from \SU(2)-structures]{Balanced Hermitian metrics from \SU(2)-structures}
\author[M. Fern\'andez, A. Tomassini, L. Ugarte and R. Villacampa]{M. Fern\'andez, A. Tomassini, L. Ugarte and R. Villacampa}
\date{\today}
\address{\noindent{\sf M. Fern\'andez:} Departamento de Matem\'aticas\\
Facultad de Ciencia  Tecnolog\'{\i}a\\ Universidad del Pa\'{\i}s Vasco\\
Apartado 644\\ 48080 Bilbao\\ Spain} \email{marisa.fernandez@ehu.es}
\address{\noindent{\sf A. Tomassini:} Dipartimento di Matematica\\ Universit\`a di Parma\\ Viale G.P. Usberti 53/A\\
43100 Parma\\ Italy} \email{adriano.tomassini@unipr.it}
\address{\noindent{\sf L. Ugarte and R. Villacampa:} Departamento de
Matem\'aticas \!-\! I.U.M.A.\\ Universidad de Zaragoza\\ Campus Plaza
San Francisco\\ 50009 Zaragoza\\ Spain}
\email{ugarte@unizar.es, raquelvg@unizar.es}
\begin{abstract}
We study the intrinsic geometrical structure of hypersurfaces
in $6$-manifolds carrying a balanced Hermitian $\SU(3)$-structure, which
we call \emph{balanced} $\SU(2)$-\emph{structures}.
We provide
conditions which imply that such a $5$-manifold can be
isometrically embedded as a hypersurface in a
manifold with a balanced $\SU(3)$-structure.
We show that any $5$-dimensional compact nilmanifold
has an invariant balanced $\SU(2)$-structure as well as new examples
of balanced Hermitian $\SU(3)$-metrics
constructed from balanced $\SU(2)$-structures.
Moreover, for $n=3,4$, we present examples of
compact manifolds, endowed with a balanced $\SU(n)$-structure,
such that the corresponding Bismut connection has holonomy equal to $\SU(n)$.

\end{abstract}
\maketitle
\section{Introduction}
Let $(J,g)$ be a Hermitian structure on a manifold $M$, with
K\"ahler form $F$ and Lie form $\theta$. The $3$-form $JdF$ is
the torsion of the Bismut connection of $(J,g)$, that is, the unique Hermitian
connection with totally skew-symmetric torsion~\cite{Bi}.
If $JdF$ is closed and non-zero (so $(M,J,g)$ is not a K\"ahler manifold)
the Hermitian structure is called {\em strong K\"ahler with torsion}.

Hermitian structures for which the Lie form $\theta$ vanishes identically or,
equivalently, $F^2 = F\wedge F$ is closed,
are called {\em balanced}.  If in addition, there is
an $\SU(3)$-structure with K\"ahler form $F$
and a complex volume $(3,0)$-form $\Psi =
\Psi_{+} + i \Psi_{-}$ on $M$ such that
$dF^2=d\Psi_+=d\Psi_-=0$, we say that
$(F, \Psi_{+}, \Psi_{-})$ is a {\em balanced $\SU(3)$-structure}.
Such structures are a generalization of integrable $\SU(3)$-structures,
which are defined by a triplet $(F, \Psi_{+}, \Psi_{-})$
satifying $dF=d\Psi_+=d\Psi_-=0$.
Balanced $\SU(3)$-structures have been very useful in physics
for the construction of explicit compact supersymmetric valid solutions
to the heterotic equations of motion in dimension six~\cite{FIUV}.

Recently, Conti and Salamon~\cite{CS} introduced the notion of
 \emph{hypo structures} on $5$-manifolds as those structures
 corresponding to the restriction of an integrable $\SU(3)$-structure
 on a $6$-manifold $M$ to a hypersurface $N$ of $M$.
They are a  generalization in dimension $5$ of
Sasakian-Einstein metrics. In fact, Sasakian-Einstein metrics
correspond to Killing spinors and {\em hypo structures} are induced
by {\em generalized} Killing spinors. In terms of differential
forms, a {\em hypo structure} on a $5$-manifold $N$ is determined by
a quadruplet $(\eta,\omega_i, 1\leq i \leq 3)$ of differential
forms, where $\eta$ is a nowhere vanishing $1$-form and $\omega_i$
are $2$-forms on $N$ satisfying certain relations (see \eqref{rhypo}
in Section~\ref{balanced 5-manifolds}).
If the forms $\eta$ and $\omega_i$ satisfy
 \begin{equation}\label{e-s}
  d\eta=-2\omega_3,  \quad d\omega_1=3\eta\wedge\omega_2,   \quad
  d\omega_2=-3\eta\wedge\omega_1,
  \end{equation}
then $N$ is a Sasakian-Einstein manifold, that is, a Riemannian
manifold such that $N\times\mathbb R$ with the cone metric is
K\"ahler and Ricci flat \cite{BGal}, so $N\times\mathbb R$ has
holonomy contained in $\SU(3)$ or, equivalently, it has an integrable
$\SU(3)$-structure.

In the general case of a hypo structure, in
\cite{CS} it is proved that a real analytic hypo structure on $N$
can be lifted to an integrable $\SU(3)$-structure on $N\times I$, for
some open interval $I$, if $(\eta,\omega_i, 1\leq i \leq 3)$ belongs
to a one-parameter family of hypo structures $(\eta(t),\omega_i(t),
1\leq i \leq 3)$ satisfying the evolution equations~\eqref{evolCS}
given in Section \ref{evol-Sect}. Moreover,  any oriented hypersurface of
a $6$-manifold with an integrable $\SU(3)$-structure is naturally
endowed with a hypo structure (see Section \ref{balanced 5-manifolds} for
details).

Our purpose in this paper is to
study the geometrical structure of an oriented hypersurface
in a $6$-manifold equipped with a balanced $\SU(3)$-structure, which
we call {\em balanced $\SU(2)$-structure}. These
are $\SU(2)$-structures defined by
a quadruplet $(\eta,\omega_i, 1\leq i \leq 3)$ of differential
forms satisfying the relations given by~(\ref{SU(2)-balanced-conditions}),
and they can be considered as a generalization to dimension five of the holomorphic
symplectic structures in dimension four. We prove
that any compact nilmanifold of dimension $5$
has a balanced $\SU(2)$-structure although three
of the $5$-nilmanifolds do not admit an invariant hypo structure~\cite{CS}. Furthermore,
we provide conditions under which a balanced $\SU(2)$-structure
on $N$ can be lifted to a balanced  $\SU(3)$-structure on $N \times\mathbb R$.
This allows us to exhibit examples of (non-invariant) balanced Hermitian metrics and,
in particular of
complex structures,
on nilmanifolds not admitting invariant balanced Hermitian metrics.

To this end, in Section~\ref{balanced 5-manifolds}
we show that there exists a balanced $\SU(2)$-structure on any oriented
hypersurface $f\colon N\longrightarrow M$ of a $6$-manifold $M$ with
a balanced $\SU(3)$-structure. Furthermore,
we prove that any balanced $\SU(2)$-structure
on $N$ can be lifted to a balanced Hermitian $\SU(3)$-structures
on $N\times\mathbb R$ if and only if it satisfies
the {\em evolution balanced $\SU(2)$ equations} (\ref{evolution-balanced})
established in Theorem~\ref{evolution-eqs} of Section~\ref{evol-Sect}.

On the other hand,
in Proposition~\ref{from-holomorphic-symplectic-to-balanced}
we characterize the circle bundles $S^1\subset N \to X$ with a balanced $\SU(2)$-structure
induced by a holomorphic
symplectic structure on $X$. In particular, we get that the
product $X\times\mathbb R$ has a balanced $\SU(2)$-structure.
Then, in Section~\ref{evol-Sect} we solve the
evolution balanced $\SU(2)$ equations (\ref{evolution-balanced})
for the compact $5$-nilmanifolds not having a hypo structure
as well as for the balanced $\SU(2)$-structures
on manifolds which are the total space of a circle bundle over Kodaira-Thurston manifold
with a certain holomorphic symplectic structure. In this way, we get
new examples of balanced Hermitian $\SU(3)$-metrics.

Finally, in Section ~\ref{Bismut} we describe in detail several examples of
compact manifolds obtained as quotient of solvable Lie groups, and
endowed with a balanced $\SU(n)$-structure,
for $n=3, 4$. In particular, we show that the holonomy of the Bismut connection
is equal to $\SU(n)$. Moreover, many of such examples are holomorphic parallelizable.

\section{Balanced \SU(2)-structures}\label{balanced 5-manifolds}
In this section we
introduce a special type of $\SU(2)$-structures on $5$-manifolds,
namely \emph{balanced $\SU(2)$-structures}, which permit construct
new Hermitian balanced metrics by
solving suitable evolution equations.
First we recall some facts about $\SU(2)$-structure on a
5-dimensional manifold. (For more details, we refer e.g. to
\cite{CS}.) Let $N$ be a 5-dimensional manifold and let $L(N)$ be
the principal bundle of linear frames on $N$. An {\em
$\SU(2)$-structure} on $N$ is an $\SU(2)$-reduction of $L(N)$. We
have the following (see \cite[Prop.1]{CS})
\begin{proposition}\label{su2structure}
$\SU(2)$-structures on a $5$-manifold $N$ are in $1:1$
correspondence with quadruplets $(\eta,\omega_1,\omega_2,\omega_3)$,
where $\eta$ is a $1$-form and $\omega_i$ are $2$-forms on $N$
satisfying
$$
  \omega_i\wedge\omega_j=\delta_{ij}v, \quad
  v\wedge\eta\not=0,
$$
for some $4$-form $v$, and
$$
i_X\omega_3=i_Y\omega_1\Rightarrow \omega_2(X,Y)\ge 0,
$$
where $i_X$ denotes the contraction by $X$.\newline Equivalently, an
$\SU(2)$-structure on $N$ can be viewed as the datum of
$(\eta,\omega_3,\Phi)$, where $\eta$ is a $1$-form, $\omega_3$ is a
$2$-form and $\Phi =\omega_1+i\omega_2$ is a complex $2$-form such
that
$$
\begin{array}{ll}
\eta\wedge\omega_3^2 \neq 0\,, &\quad \Phi^2=0\\[6pt]
\omega_3\wedge\Phi =0\,, &\quad \Phi\wedge\overline{\Phi}
=2\omega_3^2
\end{array}
$$
and $\Phi$ is of type $(2,0)$ with respect to $\omega_3$.
\end{proposition}
As a corollary of the last Proposition, we obtain the useful local
characterization of $\SU(2)$- structures (see \cite{CS}):
\begin{corollary}
If $(\eta,\omega_1,\omega_2,\omega_3)$ is an $\SU(2)$-structure on a
$5$-dimensional manifold $N$, then locally, there exists a basis of
$1$-forms $\{ e^1,\ldots ,e^5\}$ such that
$$
\eta = e^{1},\qquad \omega_1 = e^{24}+e^{53},\qquad \omega_2 =
e^{25}+e^{34},\qquad \omega_3 = e^{23}+e^{45}\,.
$$
\end{corollary}
As a consequence, $\SU(2)$-structures naturally arise on
hypersurfaces of $6$-manifolds with an $\SU(3)$-structure. Indeed,
let $f:N\longrightarrow M$ be an oriented hypersurface of a
$6$-manifold $M$ endowed with an $\SU(3)$-structure
$(F,\Psi_+,\Psi_-)$ and denote by $\mathbb U$ the unit normal vector
field. Then $N$ inherits an $\SU(2)$-structure
$(\eta,\omega_1,\omega_2,\omega_3)$ given by
\begin{equation}\label{induced-SU(2)-structure}
\eta=-i_{\mathbb U}F,\quad \omega_1=i_{\mathbb U} \Psi_-, \quad
 \omega_2=-i_{\mathbb U} \Psi_+, \quad
 \omega_3= f^*F. \quad
\end{equation}
Conversely, an $\SU(2)$-structure
$(\eta,\omega_1,\omega_2,\omega_3)$ on $N$ induces an
$\SU(3)$-structure $(F,\Psi_+,\Psi_-)$ on $N\times\mathbb R$ given
by
\begin{equation}\label{SU(3)-structure}
F=\omega_3+\eta\wedge dt, \quad\quad \Psi = \Psi_+ + i \Psi_- =
(\omega_1 + i \omega_2)\wedge(\eta+ i dt),
\end{equation}
where $t$ is a coordinate on $\mathbb R$.
\begin{definition}\label{SU(2)-balanced}
An $\SU(2)$-structure $(\eta,\omega_1,\omega_2,\omega_3)$ on a
$5$-dimensional manifold $N$ is called \emph{balanced} if it
satisfies the following equations

\begin{equation}\label{SU(2)-balanced-conditions}
d(\omega_1\wedge\eta)=0, \qquad d(\omega_2\wedge\eta)=0, \qquad
d(\omega_3\wedge\omega_3)=0.
\end{equation}
\end{definition}

\medskip

In ~\cite{CS}, an $\SU(2)$-structure is said to be {\em hypo} if
  \begin{equation}\label{rhypo}
d(\omega_1\wedge\eta)=d(\omega_2\wedge\eta)=d\omega_3=0.
 \end{equation}

\smallskip

Hence, it is obvious that any hypo structure is balanced. However,
there are nilmanifolds admitting no invariant hypo structure, but
having balanced $\SU(2)$-structures. In fact, if the Lie algebra
underlying a 5-dimensional compact nilmanifold $N$ is isomorphic to
$(0,0,0,12,14)$, $(0,0,12,13,23)$ or $(0,0,12,13,14+23)$, then there
is no invariant hypo structure on $N$~(see \cite{CS}). It is easy to check
that the $\SU(2)$-structure given by
$$
\eta = e^{1},\qquad \omega_1 = e^{24}+e^{53},\qquad \omega_2 =
e^{25}+e^{34},\qquad \omega_3 = e^{23}+e^{45},
$$
satisfies~(\ref{SU(2)-balanced-conditions}) on each one of these
three Lie algebras. Therefore, we get the following
\begin{proposition}\label{clasif}
Any 5-dimensional compact nilmanifold has an invariant
balanced $\SU(2)$-struc\-ture.
\end{proposition}
There exist also $5$-dimensional solvable non-nilpotent Lie algebras
with no invariant hypo structure, but having a balanced
$\SU(2)$-structure. For example, let us consider the Lie algebra
$\frg$ whose dual is spanned by $(e^1,\ldots ,e^5)$ such that
$$
de^1=0\,,\,\,de^2=0\,,\,\,de^3=e^{13}\,,\,\,de^4=-e^{14}\,,\,\,de^5=e^{34}\,.
$$
Then $\frg$ is a $5$-dimensional solvable non-nilpotent Lie algebra.
The simply-connected Lie group $G$ associated with $\frg$ has a
uniform discrete subgroup $\Gamma$, so that $N=\Gamma\backslash G$
is a compact solvmanifold. In fact, the manifold $N$ is the
topological product of the unit circle by the compact solvable
$4$-dimensional manifold, studied in~\cite{ACFM},
which is a circle bundle over the compact solvmanifold $Sol(3)$. A
straightforward calculation shows that the following forms
$$
\eta = e^{1},\qquad \omega_1 = e^{24}+e^{53},\qquad \omega_2 =
e^{25}+e^{34},\qquad \omega_3 = e^{23}+e^{45}\,.
$$
satisfy
$$
d(\omega_1\wedge\eta)=d(\omega_3\wedge\eta)=d(\omega_2^2)=0,
$$
and thus they define a balanced $\SU(2)$-structure on $N$.
However, $N$ has not invariant hypo $\SU(2)$-structures. First,
using Hattori's theorem~\cite{Hat}, we have that the real cohomology
groups of $N$ of degree
$\leq 2$ are
\begin{eqnarray*}
H^0(N) = \la 1\ra,\quad H^1(N) = \la [e^1], [e^2]\ra, \quad H^2(N) =
\la[e^{12}]\ra.
\end{eqnarray*}
Now, let us suppose that $N$ has an invariant hypo
$\SU(2)$-structure $(\eta,\omega_1,\omega_2,\omega_3)$. Then
$$
\omega_3 = a e^{12} + b e^{13} + c e^{14} +f e^{34},
$$
for some real numbers $a,b,c$ and $f$. Therefore, $\omega_3^2 = 2af
e^{1234}$, and so $\eta=e^5 + \sum_{i=1}^{4} {\lambda_i e^i}$ since
$\omega_3^2 \wedge \eta$ is a volume form. On the other hand,
$$
\omega_1 = \sum_{i,j=1}^{4} {a_{ij} e^{ij}},
$$
and
$$
\omega_2 = \sum_{i,j=1}^{4} {b_{ij} e^{ij}},
$$
for some real numbers $a_{ij}$ and $b_{ij}$. Now, the conditions
$d(\omega_1 \wedge \eta) = d(\omega_2 \wedge \eta) =0$ imply that
$$
\omega_1 = a_{13} e^{13}+ a_{14} e^{14}+ a_{34} e^{34},
$$
and
$$
\omega_2 = b_{13} e^{13}+ b_{14} e^{14}+ b_{34} e^{34},
$$
which implies that $\omega_1^2=\omega_2^2=0$. This is not possible
for an $\SU(2)$-structure on $N$.

\medskip

Balanced $\SU(2)$-structures on 5-manifolds are related to Hermitian
balanced structures in six dimensions as the next proposition shows.
First, recall that a balanced $\SU(3)$-structure on a 6-manifold $M$
is an $\SU(3)$-structure $(F,\Psi=\Psi_+ + i \Psi_-)$ (where $F$ is
the K\"ahler form of an almost Hermitian structure and $\Psi=\Psi_+ +
i \Psi_-$ is a complex volume form) such that $F^2$ and $\Psi$ are
closed. The latter condition implies that the underlying almost
complex structure is integrable. Notice that any balanced
$\SU(3)$-structure is in particular half-flat~\cite{ChS, H, CT}, that is,
it satisfies $dF^2 = d\Psi_+ =0$.

\begin{proposition}\label{hypersurface}
Let $f\colon N\longrightarrow M$ be an immersion of an oriented
5-manifold into a 6-manifold with an $\SU(3)$-structure. If the
$\SU(3)$-structure is balanced then
the $\SU(2)$-structure on $N$ given by
(\ref{induced-SU(2)-structure}) is
 balanced.
\end{proposition}

\begin{proof}
From (\ref{induced-SU(2)-structure}) it follows that
$\omega_1\wedge\eta= f^* \Psi_+$, $\omega_2\wedge\eta= f^* \Psi_-$
and $\omega_3\wedge\omega_3=f^* F^2$. Now, if $F^2$ and $\Psi$ are
closed then the induced structure is half-balanced.
\end{proof}

\smallskip

Let $(X,J)$ be a complex surface. By definition, a {\em holomorphic
symplectic structure} on $X$ is the datum of a $d$-closed and
non-degenerate $(2,0)$-form $\omega$ on $X$. Let $g$ be a
$J$-Hermitian metric on $X$ and $\omega_3$ be the fundamental form
of $(g,J)$. Then, up to a conformal change, we may assume that
$$
\omega_3^2=\omega_1^2=\omega_2^2\,.
$$
Then we have the following

\begin{proposition}\label{from-holomorphic-symplectic-to-balanced}
Let $(X,J)$ be a complex surface equipped with a holomorphic
symplectic structure $\omega=\omega_1+i\omega_2$, and let $\omega_3$
be the K\"ahler form of a $J$-Hermitian metric.
Then, for any integral closed 2-form $\Omega$ on $X$ annihilating
$\cos \theta\, \omega_1+ \sin\theta\, \omega_2$ and $\sin
\theta\,\omega_1 - \cos\theta\, \omega_2$ for some $\theta$, there
is a principal circle bundle $\pi\colon N\longrightarrow X$ with
connection form $\rho$ such that $\Omega$ is the curvature of $\rho$
and such that the $\SU(2)$-structure $(\eta,\omega_1^\theta,
\omega_2^\theta, \omega_3^\theta)$ on $N$ given by
$$
\begin{array}{l}\label{hypo-to-half-flat}
\eta=\rho, \\[5pt]
\omega_1^\theta= \pi^*(\cos \theta\,\omega_1+ \sin\theta\, \omega_2), \\[5pt]
\omega_2^\theta= \pi^*(-\sin\theta\,\omega_1+\cos \theta\,\omega_2), \\[5pt]
\omega_3^\theta=\pi^*(\omega_3) \\[5pt]
\end{array}
$$
is a balanced $\SU(2)$-structure.
\end{proposition}

\begin{proof}
As previously remarked, we may assume that $\omega_1^2=\omega_2^2=\omega_3^2$. Since
$d\omega_1 = d\omega_2 =0$ and $d\eta=\pi^*(\Omega)$, a simple
calculation shows that
$$
d(\omega_1^\theta\wedge \eta)= \omega_1^\theta\wedge d\eta =
\pi^*((\cos\theta\,\omega_1+ \sin\theta\, \omega_2)\wedge \Omega)=0,
$$
and
$$
d(\omega_2^\theta\wedge \eta)= \omega_2^\theta\wedge d\eta =
\pi^*((-\sin\theta\,\omega_1+ \cos\theta\, \omega_2)\wedge
\Omega)=0.
$$
The existence of a principal circle bundle in the conditions above
follows from a well known result by Kobayashi~\cite{Kobayashi}.
\end{proof}

\begin{remark}
{\rm Notice that $\Omega=0$ satisfies the hypothesis in the previous
proposition for each $\theta$ and one gets the trivial circle bundle
$N=X\times \mathbb R$ with the balanced $\SU(2)$-structure which is
the natural extension to $N$ of the holomorphic symplectic structure
on $X$. }
\end{remark}

\begin{remark}\label{geiges}
{\rm Following \cite[Def.1.1]{G}, a {\em symplectic couple} on an
oriented $4$-manifold $X$ is a pair of symplectic forms
$(\omega_1,\omega_2)$ such that $\omega_1\wedge\omega_2=0$ and
$\omega_1^2, \omega_2^2$ are volume forms defining the positive
orientation. A symplectic couple is called {\em conformal} if
$\omega_1^2=\omega_2^2$. By \cite[Thm.1.3]{G} it follows that $X$
admits a conformal symplectic couple if and only if $X$ is
diffeomorphic to ({\it a}) a complex torus, ({\it b}) a $K3$ surface
or ({\it c}) a primary Kodaira surface. In the hypothesis of prop.
\ref{hypo-to-half-flat}, $(\omega_1,\omega_2)$ defines a conformal
symplectic couple on the $4$-manifold $X$ and $\omega_3$ is a
non-degenerate $2$-form on $X$ such that
$\omega_i\wedge\omega_3=0,\,i=1,2$ and
$\omega_1^2=\omega_2^2=\omega_3^2$ }
\end{remark}

{\rm Next we illustrate the construction given in
Proposition~\ref{from-holomorphic-symplectic-to-balanced} by showing
principal circle bundles over holomorphic symplectic manifolds in
the cases ({\it a}) and ({\it c}) of Remark~\ref{geiges}. Let us
consider the closed 4-manifold $X=\Gamma\backslash G$, where the Lie
algebra $\mathfrak{g}$ of $G$ has the following structure equations
\begin{equation}\label{holomorphicsymplectic}
de^1 =0\,,\quad de^2 =0\,,\quad de^3 =0\,,\quad de^4 =-\epsilon\,
e^{23}\quad\ (\epsilon=0,1).
\end{equation}
Clearly $X$ is the Kodaira-Thurston manifold for $\epsilon=1$ and a
4-torus for $\epsilon=0$. Consider the complex structure $J$ on $X$
defined by the complex $(1,0)$-forms
$$
\varphi^1 =e^1+ie^4\,,\quad \varphi^2 =e^2+ie^3,
$$
so that
$$
\omega = \varphi^1\wedge\varphi^2= \left(e^{12}+e^{34} \right)
+i\left(e^{13}-e^{24}\right) =\omega_1+i\omega_2
$$
is a holomorphic symplectic structure on $X$. The metric $g$ on $X$
given by $g=\sum_{i=1}^4e^i\otimes e^i$ is $J$-Hermitian and the
fundamental form of $(g,J)$ is precisely $\omega_3 = e^{14}+e^{23}$,
so we are in the conditions of
Proposition~\ref{from-holomorphic-symplectic-to-balanced}. Observe
that $\omega_3$ is closed only when $X$ is a 4-torus, i.e.
$\epsilon=0$.

Now let $\Omega$ be a closed 2-form on $X$ such that
\begin{equation}\label{anikila}
\Omega\wedge (\cos \theta\, \omega_1+ \sin\theta\, \omega_2)=0=
\Omega\wedge (\sin \theta\,\omega_1 - \cos\theta\, \omega_2),
\end{equation}
for some $\theta$. A direct calculations shows that
$$
\Omega= a(e^{12}-e^{34}) + b(e^{13}+e^{24}) +(\epsilon-1)c_1\,
e^{14} + c_2\, e^{23}
$$
satisfies (\ref{anikila}) for all $\theta$, where $a,b,c_1,c_2$ are
constant. Applying
Proposition~\ref{from-holomorphic-symplectic-to-balanced}, we have a
balanced $\SU(2)$-structure on the total space $N$, which is
non-hypo if $\epsilon=1$. In this case it is easy to see that $N$ is
a compact 5-nilmanifold with underlying Lie algebra isomorphic to
either $(0,0,0,0,12)$ or $(0,0,0,12,13+24)$. }

\smallskip
%
%
%
%
%
%

\section{Evolution equations and Hermitian balanced metrics}\label{evol-Sect}
Next we establish evolution equations that allow the construction of
new balanced Hermitian metrics in dimension six from balanced
$\SU(2)$-structures in dimension five.

\begin{theorem}\label{evolution-eqs}
Let $(\eta(t),\omega_1(t),\omega_2(t),\omega_3(t))$ be a family of
$\SU(2)$-structures on a 5-manifold $N$, for $t\in I=(a,b)$. Then,
the $\SU(3)$-structure on $M=N\times I$ given by
\begin{equation}\label{SU(3)-on-NxI}
F=\omega_3(t)+\eta(t)\wedge dt, \quad\quad \Psi = (\omega_1(t) + i
\omega_2(t))\wedge(\eta(t)+ i dt),
\end{equation}
is balanced Hermitian if and only if
$(\eta(t),\omega_1(t),\omega_2(t),\omega_3(t))$ is a balanced
$\SU(2)$-structure for any $t$ in the open interval $I$, and the
following evolution equations
\begin{equation}\label{evolution-balanced}
\left\{\begin{array}{l}
\partial_t(\omega_1\wedge\eta)=-d\omega_2 ,\\[5pt]
\partial_t(\omega_2\wedge\eta)=d\omega_1 ,\\[5pt]
\partial_t(\omega_3\wedge\omega_3)=-2\,d(\omega_3\wedge\eta) ,
\end{array} \right.
\end{equation}
are satisfied.
\end{theorem}

\begin{proof}
A direct calculation shows that the $\SU(3)$-structure given
by~(\ref{SU(3)-on-NxI}) satisfies
$$
dF^2= d(\omega_3\wedge\omega_3) +
(\partial_t(\omega_3\wedge\omega_3)+2\, d(\omega_3\wedge\eta))\wedge
dt,
$$
and
$$
d\Psi= d(\omega_1\wedge\eta) - (\partial_t(\omega_1\wedge\eta)+
d\omega_2)\wedge dt + i\, d(\omega_2\wedge\eta) - i\,
(\partial_t(\omega_2\wedge\eta)- d\omega_1)\wedge dt.
$$
The forms $F^2$ and $\Psi$ are both closed if and only if
$(\eta(t),\omega_1(t),\omega_2(t),\omega_3(t))$ is a balanced
$\SU(2)$-structure for any $t\in I$, and satisfies the
equations~(\ref{evolution-balanced}).
\end{proof}

\begin{remark}\label{relation-to-CS-equations}
{\rm Notice that in the special case when the balanced
$\SU(2)$-structures in the family are hypo, the last equation
in~(\ref{evolution-balanced}) reduces to
$$
\omega_3\wedge (\partial_t\omega_3 + d\eta) =0.
$$
Hence, any solution of the hypo evolution equations introduced
in~\cite{CS}, namely
\begin{equation}\label{evolCS}
\partial_t(\omega_1\wedge\eta)=-d\omega_2 ,\quad
\partial_t(\omega_2\wedge\eta)=d\omega_1 ,\quad
\partial_t\omega_3 =- d\eta,
\end{equation}
is trivially a solution of~(\ref{evolution-balanced}), and the
resulting $\SU(3)$-structure is integrable because $F$
in~(\ref{SU(3)-on-NxI}) is closed.}
\end{remark}

Next we give some explicit solutions to the balanced evolution
equations~(\ref{evolution-balanced}). Our first example is the
5-manifold $N=X\times \mathbb{R}$, where $X$ is the Kodaira-Thurston
manifold. Let us consider the non-hypo balanced $\SU(2)$-structure
given by~\eqref{holomorphicsymplectic} for $\epsilon=1$ and with zero curvature,
namely
$$
\eta = e^{5},\qquad \omega_1 = e^{12}+e^{34},\qquad \omega_2 =
e^{13}-e^{24},\qquad \omega_3 = e^{14}+e^{23}.
$$
The family of non-hypo balanced $\SU(2)$-structures on $N$ given by
$\eta(t) = e^5$ and
$$
\omega_1(t) = e^{12}+e^{3}(e^4-t\, e^5),\ \omega_2(t) = e^{13}+
(e^4-t\, e^5)e^{2},\ \omega_3(t) = e^1(e^4-t\, e^5) + e^{23},
$$
coincides with the previous one for $t=0$ and satisfies the
evolution equations~(\ref{evolution-balanced}) for any $t\in
\mathbb{R}$.

Denote by $G$ the simply-connected nilpotent Lie group with Lie
algebra $(0,0,0,0,12)$, so that $N=\Gamma\backslash G$ for some
lattice $\Gamma$ of maximal rank in $G$. It follows from
Proposition~\ref{evolution-eqs} that the $\SU(3)$-structure on
$G\times \mathbb{R}$ given by
\begin{equation}\label{solution-(0,0,0,0,12)}
\begin{array}{lll}
F & = & e^{14} + e^{23} - t\, e^{15} +  e^5 \wedge dt,\\[6pt]
\Psi_+ & = &  e^{125} + e^{345} - \left( e^{13} - e^{24} + t\, e^{25} \right) \wedge dt ,\\[6pt]
\Psi_- & = &  e^{135} - e^{245} + \left( e^{12} + e^{34} - t\,
e^{35} \right) \wedge dt ,
\end{array}
\end{equation}
is balanced Hermitian.

\medskip

Next we find explicit solutions of~(\ref{evolution-balanced}) for
the non-hypo nilpotent Lie algebras $(0,0,0,12,14)$ and
$(0,0,12,13,23)$.

\medskip

\noindent {\bf Lie algebra (0,0,0,12,14)}: It is straightforward to
check that the family of balanced $\SU(2)$-structures given by
$$
\begin{array}{l}
\eta(t) = \sqrt[3]{{2-3t\over 2}} e^1,\\[7pt]
\omega_1(t) = {1\over 2} \left(\sqrt[3]{{2\over 2-3t}}-{2-3t\over
2}\right) e^{23} +
\sqrt[3]{{2-3t\over 2}} e^{24} - \sqrt[3]{{2\over 2-3t}} e^{35},\\[7pt]
\omega_2(t) = \sqrt[3]{{2\over 2-3t}} e^{25} + \sqrt[3]{{2-3t\over
2}} e^{34},\\[7pt]
\omega_3(t) = e^{23} -{1\over 2} \left(1- {2-3t \over
2}\sqrt[3]{{2-3t\over 2}}\right) e^{24} + e^{45},
\end{array}
$$
satisfies the evolution equations~(\ref{evolution-balanced}) for
$t\in \mathbb{R}-\{2/3\}$. Observe that the volume form of the
associated Riemannian metric along the family is given by
$$
\omega_1(t)\wedge\omega_1(t)\wedge\eta(t) = 2 \sqrt[3]{{2-3t\over
2}} e^{12345},
$$
so the orientation for $t\in (-\infty,2/3)$ is opposite to the
orientation for $t\in (2/3,\infty)$.

Let $I=(-\infty,2/3)$ and denote by $G$ the simply-connected
nilpotent Lie group with Lie algebra $(0,0,0,12,14)$. The basis of
1-forms on the product manifold $G\times I$ given by
$$
\begin{array}{llll}
\alpha^1=e^2,& \  \alpha^2=e^3,& \  \alpha^3= \sqrt[3]{{2-3t\over
2}} e^4,& \  \alpha^4= {1\over 2} \sqrt[3]{{2\over 2-3t}} (e^2 + 2
e^5) -
{2-3t\over 4} e^2,\\[5pt]
\alpha^5= \sqrt[3]{{2-3t\over 2}} e^1,& \ \alpha^6=dt,&{}&{}
\end{array}
$$
is orthonormal with respect to the Riemannian metric associated to
the balanced $\SU(3)$-structure on $G\times I$. The Hermitian
balanced structure on $G\times I$ is given by
\begin{equation}\label{solution-(0,0,0,12,14)}
\begin{array}{lll}
F & = & e^{23} - {1\over 2} e^{24} + e^{45}
+ {2-3t\over 4} \sqrt[3]{{2-3t\over 2}} e^{24} + \sqrt[3]{{2-3t\over 2}} e^1 \wedge dt,\\[10pt]
\Psi_+ & = & {1\over 2} e^{123} -e^{135} - {2-3t\over 4}
\sqrt[3]{{2-3t\over 2}} e^{123} + \sqrt[3]{{(2-3t)^2\over 4}}
e^{124} \\[7pt]
&{} & - \left( \sqrt[3]{{2\over 2-3t}} e^{25} +
\sqrt[3]{{2-3t\over 2}} e^{34} \right) \wedge dt ,\\[10pt]
\Psi_- & = &  e^{125} + \sqrt[3]{{(2-3t)^2\over 4}} e^{134} \\[7pt]
& {} & + \left( {1\over 2} \sqrt[3]{{2\over 2-3t}} e^{23} -
{2-3t\over 4} e^{23} + \sqrt[3]{{2-3t\over 2}} e^{24} -
\sqrt[3]{{2\over 2-3t}} e^{35} \right) \wedge dt.
\end{array}
\end{equation}
\medskip

\noindent {\bf Lie algebra (0,0,12,13,23)}: A direct calculation
shows that the family of balanced $\SU(2)$-structures given by
$$
\begin{array}{lll}
\eta(t) &=& {2\over 2-t} e^3,\\[6pt]
\omega_1(t) &=& {2-t\over 2} (e^{15}+e^{42}),\\[6pt]
\omega_2(t) &=& {t(2-t)(t-4)\over 4} e^{12} + {2-t\over
2}(e^{14}+e^{25}),\\[6pt]
\omega_3(t) &=& e^{12} - {t(2-t)^2(t-4)\over 8} e^{25}
-{(2-t)^2\over 4} e^{45},
\end{array}
$$
satisfies the evolution equations~(\ref{evolution-balanced}) for
$t\in \mathbb{R}-\{2\}$. Observe that the volume form of the
associated Riemannian metric along the family is given by
$$
\omega_1(t)\wedge\omega_1(t)\wedge\eta(t) = -2 e^{12345},
$$
so it remains constant.

Let $I=(-\infty,2)$ and denote by $G$ the simply-connected nilpotent
Lie group with Lie algebra $(0,0,12,13,23)$. The basis of 1-forms on
the product manifold $G\times I$ given by
$$
\begin{array}{llll}
\alpha^1=e^1, &\quad \alpha^2=e^2, &\quad \alpha^3={2-t\over 2}
e^5,&
\quad \alpha^4= {t(2-t)(t-4)\over 4} e^2 + {2-t\over 2} e^4,\\[6pt]
\alpha^5= {2\over 2-t} e^3,& \quad \alpha^6=dt,&{}&{}
\end{array}
$$
is orthonormal with respect to the Riemannian metric associated to
the balanced $\SU(3)$-structure on $G\times I$. The Hermitian
balanced structure on $G\times I$ is given by
\begin{equation}\label{solution-(0,0,12,13,23)}
\begin{array}{lll}
F & = & e^{12} - {t(2-t)^2(t-4)\over 8} e^{25} - {(2-t)^2\over 4}
e^{45}
+ {2\over 2-t} e^3 \wedge dt,\\[8pt]
\Psi_+ & = & -e^{135}+e^{234} - {2-t\over 2}({t(t-4)\over 2}e^{12}
+e^{14} +e^{25})\wedge dt,\\[8pt]
\Psi_- & = & -e^{134}-e^{235} + {t(t-4)\over 2}e^{123} +  {2-t\over
2}(e^{15}-e^{24})\wedge dt.
\end{array}
\end{equation}

\smallskip

As a consequence of our previous examples we conclude:

\begin{corollary}\label{balanced Lie-groups}
The $6$-dimensional simply-connected nilpotent Lie groups $H_8$,
$H_{16}$ and $H_{17}$ corresponding to the Lie algebras
$\frh_8=(0,0,0,0,0,12)$, $\frh_{16}=(0,0,0,12,14,24)$ and
$\frh_{17}=(0,0,0,0,12,15)$ have balanced Hermitian
$\SU(3)$-structures.
\end{corollary}

It is worthy to remark that $H_8$ and $H_{16}$ have invariant
complex structures, but none of them admit invariant compatible
balanced metric~\cite{U}. On the other hand, $H_{17}$ has no
invariant complex structures.
%
%

\smallskip
%
%
%

\medskip

Finally, recall that according to \cite{H} if $M$ is a $6$-manifold
with a family $(F(s),\Psi_+(s),\Psi_-(s))$ of half-flat structures,
$s \in I=(a,b)$, satisfying the evolution equations
$$
\partial_s \Psi_+ =  d F, \quad\quad F \wedge \partial_s F = -  d \Psi_-
,
$$
then the product manifold $M\times I$ has a Riemannian metric whose
holonomy is contained in $G_2$; in fact, the form $\varphi=F(s)
\wedge ds + \Psi_+(s)$ defining the $G_2$-structure is parallel.
Therefore, the balanced Hermitian
$\SU(3)$-structures~(\ref{solution-(0,0,0,0,12)}),~(\ref{solution-(0,0,0,12,14)})
and~(\ref{solution-(0,0,12,13,23)}) can be lifted to a metric with
holonomy in $G_2$.

%
%
%

\section{Holonomy of Bismut connection of balanced Hermitian metrics on solvmanifolds}\label{Bismut}
In this section we provide some examples of compact Hermitian manifolds, endowed
with an $\SU(n)$-structure whose associated metric is
balanced and such that the corresponding Bismut connection has
holonomy equal
to $\SU(n)$, for $n=3$ and $4$.

\smallskip

Let $(M,J,g)$ be a Hermitian manifold and $F$ be the
K\"ahler form of the Hermitian structure $(g,J)$. Denote by $\nabla^{LC}$ the Levi Civita connection of the metric $g$. Then the Bismut connection
$\nabla^B$ of $(M,J,F,g)$ is characterized by the following formula
$$
g(\nabla^B_XY,Z) = g(\nabla^{LC}_XY,Z)+\frac{1}{2}T(X,Y,Z)\,,\quad\forall X,Y,Z\in\Gamma(M,TM)\,,
$$
where the torsion form $T$ is given by
$$
T(X,Y,Z) = JdF(X,Y,Z)\,.
$$
We will need to compute the curvature of $\nabla^B$. In order to do this, we will use the Cartan structure equations,
\begin{equation}\label{connection}
\left\{
\begin{array}{ll}
de^i + \displaystyle\sum_{j=1}^{2n}\omega^i_j\wedge e^j = \tau^i\,,&\qquad i=1,\ldots ,2n\\[15pt]
\omega^{i}_j+\omega^{j}_{2n} =0\,,&\qquad i,j=1,\ldots ,2n\,,
\end{array}
\right.
\end{equation}
\begin{equation}\label{curvature}
\left\{
\begin{array}{ll}
d\omega^i_j + \displaystyle\sum_{r=1}^{2n}\omega^i_r\wedge \omega^r_j = \Omega^i_j\,,&\qquad i,j=1,\ldots ,2n\\[15pt]
\Omega^{i}_j+\Omega^{j}_{2n} =0\,,&\qquad i,j=1,\ldots ,2n\,,
\end{array}
\right.
\end{equation}
where $\{e^1,\ldots ,e^{2n}\}$ is an orthonormal coframe,  $\omega^i_j$ are
the connection $1$-forms, $\tau^i$ are the torsion $2$-forms and $\Omega^i_j$ are the curvature
$2$-forms.\newline
If $\{e_1,\ldots ,e_{2n}\}$ denotes the dual frame of $\{e^1,\ldots ,e^{2n}\}$, then
$$
\tau^i =\sum_{j<k=1}^{2n} T_{i\,jk}e^j\wedge e^k\,,\quad i=1,\ldots ,2n\,,
$$
where $T_{i\,jk}=T(e_i,e_j,e_k)$.
\medskip

\subsection{Dimension six}
There are two 3-dimensional complex-parallelizable (non-Abelian)
solvable Lie groups~\cite{Nak}. The complex structure equations are
given by:
\begin{enumerate}
\item[{\rm (I)}] $d\varphi^1=d\varphi^2=0,\ d\varphi^3=\varphi^{12}$;\\
\item[{\rm (II)}] $d\varphi^1=0,\ d\varphi^2=\varphi^{12},\
d\varphi^3=-\varphi^{13}$.
\end{enumerate}
Here (I) is the Lie group underlying the {\em Iwasawa manifold}.
In the following we show that they can be endowed with a
balanced Hermitian $\SU(3)$-structure and determine the holonomy of the
corresponding Bismut connection.

\begin{theorem}\label{3-dim-complex-solvable}
Any 3-dimensional complex-parallelizable (non-Abelian) solvable Lie
group has a Hermitian metric such that the holonomy of its Bismut
connection is equal to $\SU(3)$.
\end{theorem}
\begin{proof}
First we show that for the standard Hermitian balanced structure $(J,g)$
on the Lie group (I) the holonomy of its Bismut connection equal to $\SU(3)$.

Let us consider the real basis $(e^1,\ldots,e^6)$ given by
$$
\varphi^1=e^1+i\,e^2,\quad \varphi^2=e^3+i\,e^4,\quad \varphi^3=e^5+i\,e^6.
$$
In this case the structure equations are
$$
d e^1= d e^2= d e^3= d e^4= 0,\quad d e^5= e^{13}-e^{24},\quad d e^6
= e^{14}+e^{23}.
$$
Then, the complex structure $J$ is given by
$$
J e^1=-e^2, \quad J e^2=e^1, \quad J e^3=-e^4, \quad J e^4=e^3,
\quad J e^5=-e^6, \quad J e^6=e^5.
$$
The fundamental form $F$ associated with the $J$-Hermitian metric $g=\sum_{i=1}^6e^i\otimes e^i$
is given by
$$
F=e^{12}+e^{34}+e^{56}.
$$
Since $dF=e^{136}-e^{145}-e^{235}-e^{246}$, we get that $g$ is
balanced and the torsion $T=JdF$ is given by
$$
T=-e^{135}-e^{146}-e^{236}+e^{245}.
$$
The non-zero connection 1-forms of the Bismut connection $\nabla^B$
are
$$
\begin{array}{llll}
\omega^1_5=-e^3,&\quad \omega^1_6=-e^4,&\quad \omega^2_5= e^4,&\quad\omega^2_6=-e^3,\\[8pt]
\omega^3_5=e^1,&\quad \omega^3_6=e^2,&\quad\omega^4_5=-e^2,&\quad \omega^4_6=e^1.
\end{array}
$$
The (linearly independent) curvature forms of the Bismut connection
are
$$
\Omega^1_2=2e^{34},\quad\ \Omega^1_3=-e^{13}-e^{24},\quad\
\Omega^2_3=e^{14}-e^{23},\quad\ \Omega^3_4=2e^{12},
$$
and the (linearly independent) curvature forms of the Bismut
connection are
$$
\begin{array}{ll}
\nabla^B_{E_1}(\Omega^1_2)=-2(e^{36}-e^{45}),&\quad\nabla^B_{E_2}(\Omega^1_2)=2(e^{35}+e^{46}),\\[8pt]
\nabla^B_{E_3}(\Omega^3_4)=2(e^{16}-e^{25}),&\quad\nabla^B_{E_4}(\Omega^3_4)=-2(e^{15}+e^{26}).
\end{array}
$$
Therefore, Hol$(\nabla^B)=\SU(3)$.

\vskip.3cm

For case (II) we consider the real basis $(e^1,\ldots,e^6)$ given by
$$
\varphi^1=e^1+i\,e^2,\quad \varphi^2=e^3+i\,e^4,\quad \varphi^3=e^5+i\,e^6.
$$
In terms of this basis, the structure equations are
$$
\begin{array}{lll}
d e^1= d e^2= 0,&\quad d e^3=e^{13}-e^{24} ,&\quad d e^4=e^{14}+e^{23},\\[6pt]
d e^5= -e^{15}+e^{26},&\quad d e^6 = -e^{16}-e^{25}.&{}
\end{array}
$$
Then, the complex structure $J$ is given by
$$
J e^1=-e^2, \quad J e^2=e^1, \quad J e^3=-e^4, \quad J e^4=e^3,
\quad J e^5=-e^6, \quad J e^6=e^5.
$$
Then the fundamental form $F$ associated with the $J$-Hermitian metric\newline
$g=\sum_{i=1}^6e^i\otimes e^i$ is given by
$$
F=e^{12}+e^{34}+e^{56}.
$$
Since $dF=2(e^{134}-e^{156})$, we get that $g$ is balanced and
the torsion $T=JdF$ is given by
$$
T=-2(e^{234}-e^{256}).
$$
The non-zero connection 1-forms of the Bismut connection $\nabla^B$
are
$$
\begin{array}{llll}
\omega^1_5=-e^3,&\quad \omega^1_6=-e^4,&\quad \omega^2_5= e^4,&\quad \omega^2_6=-e^3,\\[8pt]
\omega^3_5=e^1,&\quad \omega^3_6=e^2,&\quad\omega^4_5=-e^2,&\quad \omega^4_6=e^1.
\end{array}
$$
The (linearly independent) curvature forms of the Bismut connection
are
$$
\begin{array}{llll}
\Omega^1_3=-e^{13}-e^{24},& \Omega^1_4=-e^{14}+e^{23},& \Omega^1_5=-e^{15}-e^{26},&\Omega^1_6=-e^{16}+e^{25}\\[8pt]
\Omega^3_4=-2e^{34},& \Omega^3_5=e^{35}+e^{46},& \Omega^3_6=e^{36}-e^{45},& \Omega^5_6=-2e^{56}.
\end{array}
$$
Therefore, the holonomy of $\nabla^B$ is again equal to $\SU(3)$.
\end{proof}

\smallskip











Next, we study the holonomy of the Bismut connection
for Hermitian balanced metrics on $6$-nilmanifolds.
First of all, we recall ~\cite{U} that if a nilmanifold $M=\Gamma\backslash G$ admits a Hermitian balanced
metric (not necessarily invariant), then the Lie algebra $\frg$ of
$G$ is isomorphic to $\frh_{1},\ldots,\frh_6$ or $\frh_{19}^-$,
where $\frh_{1} = (0,0,0,0,0,0)$ is the Abelian Lie algebra and
$$
\begin{array}{rcl}
\frh_{2} &\!\!=\!\!& (0,0,0,0,12,34),\\[2pt]
\frh_{3} &\!\!=\!\!& (0,0,0,0,0,12+34),\\[2pt]
\frh_{4} &\!\!=\!\!& (0,0,0,0,12,14+23),
\end{array}
\quad\quad\quad
\begin{array}{rcl}
\frh_{5} &\!\!=\!\!& (0,0,0,0,13+42,14+23),\\[2pt]
\frh_{6} &\!\!=\!\!& (0,0,0,0,12,13),\\[2pt]
\frh^-_{19} &\!\!=\!\!& (0,0,0,12,23,14-35).
\end{array}
$$

\noindent where $\frh_5$ is the Lie algebra underlying the Iwasawa manifold
considered in Theorem~\ref{3-dim-complex-solvable}.
Therefore, the Lie algebra $\frh_{19}^-$ is the unique 6-dimensional 3-step
nilpotent Lie algebra admitting balanced structures,

\smallskip

On the other hand, in~\cite{FPS}[Corollary 6.2, Example 6.3] an explicit balanced Hermitian metric on
$\frh_{3},\frh_{4},\frh_{5},\frh_{6}$ is given such that the holonomy of the corresponding
Bismut connection is equal to $\SU(3)$ for $\frh_{4},\frh_{5}$ and $\frh_{6}$. For $\frh_3$ they
show that there is a reduction of the holonomy group, whereas it remains to study the Lie algebra
$\frh_2$ because the particular coefficients given at the end of the proof of~\cite{FPS}[Corollary 6.2]
correspond to $\frh_6$.

Next we prove that $\frh_{2}$ and $\frh_{19^{-}}$ also admit a balanced Hermitian
$\SU(3)$-structure such that the holonomy of its Bismut connection also equals $\SU(3)$.
That is, we prove the  following

\begin{theorem}
Let $M$ be a 6-dimensional nilmanifold admitting invariant balanced
Hermitian structures $(g,J)$. If the first Betti number of $M$ is
$\leq 4$, then there is $(g,J)$ such that the holonomy of its Bismut
connection is equal to $\SU(3)$.
\end{theorem}

\begin{proof}
A 6-dimensional nilmanifold with $b_1(M)\leq 4$ admitting invariant
balanced Hermitian structures $(g,J)$ has underlying Lie algebra
isomorphic to $\frh_2$, $\frh_4$, $\frh_5$, $\frh_6$ or
$\frh_{19}^-$.

In order to
define a balanced Hermitian $\SU(3)$-structure on $\frh_2$ such that
the holonomy of its Bismut connection is equal to $\SU(3)$,  let us
consider the structure equations
\begin{equation}\label{ecus-h2}
d e^1= d e^2= d e^3= d e^4= 0,\quad d e^5= e^{13}-e^{24},\quad d e^6
= -2e^{12}+e^{14}+e^{23}+2e^{34},
\end{equation}
and the complex structure $J$ given by
$$
J e^1=-e^2, \quad J e^2= e^1, \quad J e^3=-e^4, \quad J e^4=e^3,
\quad J e^5=-e^6, \quad J e^6= e^5.
$$
Firstly, we notice that the structure equations~(\ref{ecus-h2}) correspond to the Lie algebra
$\frh_2$. To see this it is sufficient to express the equations with
respect to the new basis
$$
\begin{array}{lll}
f^1=-2\,e^2+\sqrt{3}\,e^3+e^4, \ \  &f^2=e^1-\sqrt{3}\,e^2+2\,e^3,\ \ &f^3= 2\,e^2+\sqrt{3}\,e^3 - e^4,\\[6pt]
f^4= e^1+\sqrt{3}\,e^2+2\,e^3,\ \ & f^5=-\sqrt{3}\,e^5-e^6,\ \ &f^6=-\sqrt{3}\,e^5+e^6.
\end{array}
$$

Now, the fundamental form $F$ associated to the $J$-Hermitian metric $g=\sum_{i=1}^6e^i\otimes e^i$
is given by $F=e^{12}+e^{34}+e^{56}$. Since
$$
dF=e^{136}-e^{246}+2e^{125}-e^{145}-e^{235}-2e^{345}\,,
$$
we get that $g$ is balanced and the torsion $T$ is given by
$$
T=-2e^{126}-e^{135}-e^{146}-e^{236}+e^{245}+2e^{346}.
$$
The (linearly independent) curvature forms for the Bismut connection
are
$$
\begin{array}{ll}
\Omega^1_2=-2(2e^{12}-e^{14}-e^{23}-3e^{34}),\quad & \Omega^1_3=-(e^{13}+e^{24}),\\[8pt]
\Omega^1_4=-(e^{14}-e^{23}),\quad & \Omega^1_5=-2e^{46},\\[8pt]
\Omega^1_6=2e^{36},\quad & \Omega^3_4=2(3e^{12}-e^{14}-e^{23}-2e^{34}), \\[8pt]
\Omega^3_5=-2e^{26},\quad & \Omega^3_6=2e^{16}.
\end{array}
$$
Therefore, Hol$(\nabla^B)=\SU(3)$.

We have given a balanced metric with Hol$(\nabla^B)=\SU(3)$ for each case.
\end{proof}

To complete the proof,  we show
that there is a balanced Hemitian $\SU(3)$-structure $(J,g)$ on $\frh_{19}^-$
with holonomy of its Bismut connection equal to $\SU(3)$.

We know that the structure equations of $\frh_{19}^-$ are
$$
d\, e^1= d\, e^2= d\, e^3= 0,\quad d\, e^4= e^{12},\quad d\, e^5=
e^{23},\quad d\, e^6 = e^{14}-e^{35}.
$$
We define the complex structure $J$ given by
$$
J e^1=e^3, \quad J e^2=4\, e^6, \quad J e^3=-e^1, \quad J e^4=-e^5,
\quad J e^5=e^4, \quad J e^6=-\frac14\, e^2.
$$

Then, the fundamental form $F$ associated to the $J$-Hermitian metric $g=\sum_{i=1}^6e^i\otimes e^i$
is given by $F=-e^{13}-e^{26}+e^{45}$ and
$$
dF=-e^{124}+e^{125}-e^{234}-e^{235}\,.
$$
Therefore, we get that $g$ is balanced and the torsion $T$ is given by
$$
T=e^{146}-e^{156}-e^{346}-e^{356}.
$$
The non-zero connection 1-forms of the Bismut connection $\nabla^B$
are
$$
\begin{array}{llll}
\omega^1_2=-\frac12 e^4,\quad &\omega^1_4=-\frac12 e^2-e^6,\quad & \omega^1_5=\frac12 e^6,\quad & \omega^1_6=-\frac12 e^5,\\[8pt]
\omega^2_3=-\frac12 e^5,\quad & \omega^2_4=\frac12 e^1, &\omega^2_5=-\frac12 e^3,\quad & \omega^3_4=\frac12 e^6\\[8pt]
\omega^3_5=\frac12 e^2+ e^6,\quad & \omega^3_6=-\frac12 e^4,\quad & \omega^4_6=\frac12 e^3,\quad & \omega^5_6=\frac12 e^1.
\end{array}
$$
and the (linearly independent) curvature forms for the Bismut
connection are
$$
\begin{array}{ll}
\Omega^1_2=-\frac14 (3e^{12}+2e^{16}+e^{36}),\quad & \Omega^1_3=\frac12 (e^{26}+e^{45}),\\[8pt]
\Omega^1_4=-\frac34(e^{14}-e^{35}),\quad & \Omega^1_5=\frac14 (2e^{14}-e^{15}- e^{34}-2 e^{35}),\\[8pt]
\Omega^1_6=-\frac14 (e^{16}+3e^{23}+2e^{36}),\quad & \Omega^2_4=\frac14(e^{24}-2e^{46}-e^{56}),\\[8pt]
\Omega^2_5=\frac14 (e^{25}+e^{46}-2 e^{56}),\quad & \Omega^2_6=\frac12 (e^{13}-e^{45}).
\end{array}
$$
Therefore, Hol$(\nabla^B)=\SU(3)$.
\medskip\par



\smallskip

To finish this section, we show an example of a six-dimensional
compact solvmanifold with a balanced metric such that the holonomy
of its Bismut connection is equal to $\SU(3)$.\newline
Let $\mathfrak{g}$ be the solvable Lie algebra whose structure
equations are
$$
de^1=0,\;\; de^2=0,\;\; de^3 =e^1\wedge e^3,\;\; de^4 =-e^1\wedge e^4,\;\; de^5 =e^1\wedge e^5,\;\; de^6 =-e^1\wedge e^6\,.
$$
Let $G$ be the simply-connected Lie group whose Lie algebra is $\mathfrak{g}$.
It can be easily checked that, for any $X\in \mathfrak{g}$,  $\rm{ad}_X$ has
real eigenvalues, i.e. $\mathfrak{g}$ is completely solvable. Thus $G$ has
a uniform discrete
subgroup $\Gamma$, such that
$M^6=\Gamma\backslash G$ is a $6$-dimensional compact solvmanifold. \newline
Define an (integrable) complex structure $J$ on $M^6$ by setting
$$
Je^1=-e^2,\;\;Je^2=e^1,\;\;Je^3=-e^5,\;\;Je^4=-e^6,\;\;Je^5=e^3,\;\;Je^6=e^4,
$$
and a $J$-Hermitian structure $g$ on $M$ by $g =\sum_{i=1}^6 e^i\otimes e^i$. Then $F=\sum_{i=1}^3e^{2i-1}\wedge e^{2i}$. Hence
$$
dF= 2\left(e^{135}-e^{146}\right)
$$
and consequently, the torsion form of the Bismut connection is given by
$$
T= 2\left(e^{246}-e^{235}\right)\,.
$$
By the above expression, the only non-zero components of the torsion are
$$
T_{235}=-2\,,\quad T_{246}=2\,.
$$
By solving \eqref{connection}, we obtain that the non zero connection forms are given by
$$
\begin{array}{llll}
\omega^1_3=-e^3\,,&\,\,\omega^1_4=e^4\,,&\,\,\omega^1_5=-e^5\,,&\,\,\omega^1_6=e^6\,,\\[8pt]
\omega^2_3=e^5\,,&\,\,\omega^2_4=-e^6\,,&\,\,\omega^2_5=-e^3\,,&\,\,\omega^2_6=e^4\,,\\[5pt]
\omega^3_4=e^5\,,&\,\,\omega^3_5=e^2\,,&\,\,\omega^4_5=0\,,&\,\,\omega^4_6=-e^2
\end{array}
$$
Therefore, we get
$$
\begin{array}{llll}
\Omega^1_2=2\left(e^{35}+e^{46}\right),&\Omega^1_3=-e^{13}-e^{25},&\Omega^1_4=-e^{14}-e^{26},&\Omega^1_5=-e^{15}+e^{23},\\[8pt]
\Omega^1_6=-e^{16}+e^{24},&\Omega^2_3=e^{15}-e^{23},&\Omega^2_4=e^{16}-e^{24},&\Omega^2_5=-e^{13}-e^{25},\\[8pt]
\Omega^2_6=-e^{14}-e^{26},&\Omega^3_4=e^{34}+e^{56},&\Omega^3_5=-2e^{35},&\Omega^3_6=e^{36}+e^{45}\\[8pt]
\Omega^4_5=e^{36}+e^{45},&\Omega^4_6=-2e^{46},&\Omega^5_6=e^{34}+e^{56}\,.{}&{}
\end{array}
$$
By the above expression for the curvature forms, we obtain that $\Hol(\nabla^B)=\SU(3)$. On the other hand, the manifold $M$ has no
K\"ahler structures. Indeed, in view of the Main Theorem in \cite{Ha}, a
compact solvmanifold of completely solvable type has a K\"ahler
structure if and only if it is a complex torus. Finally, note that
$$
dF^2=2\left(e^{135}-e^{146}\right)\wedge\left(e^{12} + e^{35} + e^{46}\right)=0\,,
$$
i.e. $g$ is a balanced metric.

\subsection{Balanced Hermitian metrics on compact 8-solvmanifolds}
Here we present examples of compact balanced Hermitian $8$-manifolds
with an $\SU(4)$-structure such that the holonomy of its Bismut
connection is all the Lie group $\SU(4)$. Furthermore, according to \cite{Nak}, all the examples are holomorphic parallelizable.
\begin{example}{\rm
Consider the Lie algebra defined by the complex structure equations
$$d\varphi^1=d\varphi^2=0,\quad d\varphi^3=-\varphi^{12},\quad
d\varphi^4=-2\varphi^{13}.$$ Let $e^1,\dots,e^8$ be the basis given
by $e^{2j-1}+i\, e^{2j} = \varphi^j$, for $j=1,\ldots,4$.

\vskip.2cm

The corresponding real structure equations are
$$
\begin{array}{lll}
de^1=de^2=de^3=de^4=0,&\quad de^5=-e^{13}+e^{24},&{}\\[8pt]
de^6=-e^{14}-e^{23},&\quad de^7=-2(e^{15}-e^{26}),&\quad de^8=-2(e^{16}+e^{25}).
\end{array}
$$
Let $J$ be the complex structure given by
$$
\begin{array}{llll}
Je^1=-e^2,\quad & Je^2=e^1,\quad & Je^3=-e^4,\quad & Je^4=e^3,\\[6pt]
Je^5=-e^6, \quad & Je^6=e^5,\quad & Je^7=-e^8,\quad & Je^8=e^7.
\end{array}
$$

Then fundamental form $F$ associated with the $J$-Hermitian metric\newline
$g=\sum_{i=1}^8e^i\otimes e^i$
is given by $F=\sum_{j=1}^4e^{2j-1}\wedge e^{2j}$.
Since
$$
dF=-e^{136}+e^{145}-2e^{158}+2e^{167}+e^{235}+e^{246}+2e^{257}+2e^{268}\,,
$$
we get that $g$ is balanced and the torsion $T$ is given by
$$T=JdF=e^{135}+e^{146}+2e^{157}+2e^{168}+e^{236}-e^{245}+2e^{258}-2e^{267}.$$
The non-zero connection 1-forms of the Bismut connection $\nabla^B$
are
$$
\begin{array}{llllll}
\omega^1_5=e^3,& \omega^1_6=e^4,&\omega^1_7=2e^5,&\omega^1_8=2e^6,& \omega^2_5=-e^4,& \omega^2_6=e^3,\\[8pt]
\omega^2_7=-2e^6,& \omega^2_8=2e^5,& \omega^3_5=-e^1,& \omega^3_6=-e^2,& \omega^4_5=e^2,&\omega^4_6=-e^1,\\[8pt]
\omega^5_7=-2e^1,& \omega^5_8=-2e^2,& \omega^6_7=2e^2,& \omega^6_8=-2e^1.{}&{}&
\end{array}
$$
The following curvature forms for the Bismut connection are linearly
independent:
$$
\begin{array}{lll}
\Omega^1_3=-e^{13}-e^{24},\quad & \Omega^1_4=-e^{14}+e^{23},\quad & \Omega^1_5=-4(e^{15}+e^{26}),\\[8pt]
\Omega^1_6=-4(e^{16}-e^{25}),\quad & \Omega^3_4=2e^{12},\quad & \Omega^5_6=2(3e^{12}-e^{34}),\\[8pt]
\Omega^5_7=-2(e^{35}+e^{46}),\quad &\Omega^5_8=-2(e^{36}-e^{45}),\quad &\Omega^7_8=-8(e^{12}+e^{56}).
\end{array}
$$
This gives a 9-dimensional space. Moreover, the following 6
covariant derivatives of the curvature forms
$$
\begin{array}{ll}
\nabla^B_{E_2}(\Omega^1_2)=-16(e^{57}+e^{68}),\!\!&\!\! \nabla^B_{E_5}(\Omega^1_3)=-2(e^{37}+e^{48}),
\nabla^B_{E_6}(\Omega^1_3)=-2(e^{38}-e^{47}),\\[8pt]
\nabla^B_{E_6}(\Omega^1_5)=-8(e^{58}-e^{67}),\!\!&\!\! \nabla^B_{E_5}(\Omega^3_4)=-4(e^{18}-e^{27}),
\nabla^B_{E_6}(\Omega^3_4)=4(e^{17}+e^{28}),
\end{array}
$$
are linearly independent, therefore Hol($\nabla^B$)=\SU(4).}
\end{example}


\begin{example}
{\rm Consider the complex structure equations
$$
d\varphi^1=0,\quad d\varphi^2=\varphi^{12},\quad d\varphi^3=c\varphi^{13},\quad
d\varphi^4=-(1+c)\varphi^{14},
$$
where $c(1+c)\neq0$. Let
$e^1,\dots,e^8$ be the basis given by $\varphi^j=e^{2j-1}+i\, e^{2j}$, for $j=1,\ldots,4$.

\vskip.2cm\noindent
The corresponding real structure equations are
$$
\begin{array}{llll}
de^1=de^2=0,& de^3=e^{13}-e^{24},&de^4=e^{14}+e^{23},{}&\\[6pt]
de^5=c(e^{15}-e^{26})& de^6=c(e^{16}+e^{25}),& de^7=-(1+c)(e^{17}-e^{28}),{}& \\[6pt]
de^8=-(1+c)(e^{18}+e^{27})\,.{}&{}&{}&
\end{array}
$$

Let $J$ be the complex structure given by
$$
\begin{array}{llll}
Je^1=-e^2,\quad & Je^2=e^1,\quad & Je^3=-e^4,\quad & Je^4=e^3,\\[6pt]
Je^5=-e^6,\quad & Je^6=e^5,\quad & Je^7=-e^8,\quad & Je^8=e^7.
\end{array}
$$
Then $g=\sum_{i=1}^8e^i\otimes e^i$ is a $J$-Hermitian metric and the fundamental form
is given by $F=\sum_{j=1}^4e^{2j-1}\wedge e^{2j}$. Hence
$$
dF=2\left(e^{134}+ce^{156}-(1+c)e^{178}\right)
$$
and consequently we get that $g$ is
balanced and the torsion $T$ is expressed by
$$
T=JdF=-2\left(e^{234}+ce^{256}-(1+c)e^{278}\right).
$$
A direct computation shows that the non-zero connection 1-forms of the Bismut connection $\nabla^B$
are
$$
\begin{array}{lllll}
\omega^1_3=-e^3,& \omega^1_4=-e^4,&\omega^1_5=-ce^5,& \omega^1_6=-ce^6,& \omega^1_7=(1+c)e^7,\\[8pt]
\omega^1_8=(1+c)e^8,& \omega^2_3=e^4,& \omega^2_4=-e^3,& \omega^2_5=ce^6,& \omega^2_6=-ce^5,\\[8pt]
\omega^2_7=-(1+c)e^8,& \omega^2_8=(1+c)e^7,&\omega^3_4=2e^2,& \omega^5_6=2ce^2,& \omega^7_8=-2(1+c)e^2.
\end{array}
$$
The following curvature forms for the Bismut connection are linearly
independent:
$$
\begin{array}{lll}
\Omega^1_3=-(\alpha^{13}+\alpha^{24}),&\Omega^1_4=-\alpha^{14}+\alpha^{23}, &\Omega^1_5=-c^2(\alpha^{15}+\alpha^{26}),\\[8pt]
\Omega^1_6=-c^2(\alpha^{16}-\alpha^{25}), &\Omega^1_7=-(1+c)^2(\alpha^{17}+\alpha^{28}), &{}\\[8pt]
\Omega^1_8=-(1+c)^2(\alpha^{18}-\alpha^{27}), & {} &\\[8pt]
\Omega^3_4=-2\alpha^{34},& \Omega^3_5=-c(\alpha^{35}+\alpha^{46}), & \Omega^3_6=-c(\alpha^{36}-\alpha^{45}),\\[8pt]
\Omega^3_7=(1+c)(\alpha^{37}+\alpha^{48}),& \Omega^3_8=(1+c)(\alpha^{38}-\alpha^{47}),&\Omega^5_6=-2c^2\alpha^{56},\\[8pt]
\Omega^5_7=c(1+c)(\alpha^{57}+\alpha^{68}), & \Omega^5_8=c(1+c)(\alpha^{58}-\alpha^{67}), & \Omega^7_8=-2(1+c)^2\alpha^{78}.
\end{array}
$$
Therefore, Hol($\nabla^B$)=\SU(4) for any value of the parameter $c$.
}
\end{example}


\begin{example}
{\rm
Consider the complex structure equations
$$
d\varphi^1=0,\quad d\varphi^2=\varphi^{12},\quad d\varphi^3=-2\varphi^{13},\quad
d\varphi^4=\varphi^{14}-\varphi^{12}.
$$
Let $e^1,\dots,e^8$ be the
basis given by $\varphi^j=e^{2j-1}+i\, e^{2j}$, for
$j=1,\ldots,4$.

\vskip.2cm\noindent

The corresponding real structure equations are
$$
\begin{array}{ll}
de^1=de^2=0,\,de^3=e^{13}-e^{24},& de^4=e^{14}+e^{23},\,de^5=-2(e^{15}-e^{26}),\\[8pt]
de^6=-2(e^{16}+e^{25}),& de^7=-e^{13}+e^{17}+e^{24}-e^{28}, \\[8pt]
de^8=-e^{14}+e^{18}-e^{23}+e^{27}.{}&
\end{array}
$$

Let $J$ be the complex structure given by
$$
\begin{array}{llll}
Je^1=-e^2,\quad & Je^2=e^1,\quad & Je^3=-e^4,\quad & Je^4=e^3,\\[6pt]
Je^5=-e^6,\quad & Je^6=e^5,\quad & Je^7=-e^8,\quad & Je^8=e^7.
\end{array}
$$
The fundamental form $F$ associated to the $J$-Hermitian metric $g=\sum_{i=1}^8e^i\otimes e^i$ is given by
$F=\sum_{j=1}^4e^{2j-1}\wedge e^{2j}$. In such a case,
$$
dF=2e^{134}-e^{138}+e^{147}-4e^{156}+2e^{178}+e^{237}+e^{248}
$$
and the torsion $T$ is given by
$$
T=JdF=e^{137}+e^{148}-2e^{234}+e^{238}-e^{247}+4e^{256}-2e^{278}.
$$
Again, $g$ is balanced. The non-zero connection 1-forms of the Bismut connection $\nabla^B$
are
$$
\begin{array}{lllll}
\omega^1_3=-e^3,& \omega^1_4=-e^4,& \omega^1_5=2e^5,& \omega^1_6=2e^6,& \omega^1_7=e^3-e^7,\\[8pt]
\omega^1_8=e^4-e^8,& \omega^2_3=e^4,& \omega^2_4=-e^3,& \omega^2_5=-2e^6,& \omega^2_6=2e^5,\\[8pt]
\omega^2_7=-e^4+e^8,& \omega^2_8=e^3-e^7,& \omega^3_4=2e^2,&\omega^3_7=-e^1,{}&\\[8pt]
\omega^3_8=-e^2,& \omega^4_7=e^2,& \omega^4_8=-e^1,& \omega^5_6=-4e^2,& \omega^7_8=2e^2.
\end{array}
$$
The following 15 curvature forms for the Bismut connection are
linearly independent:
$$
\begin{array}{lll}
\Omega^1_3=-2e^{13}+e^{17}-2e^{24}+e^{28},& \Omega^1_4=-2e^{14}+e^{18}+2e^{23}-e^{27},&{}\\[8pt]
\Omega^1_5=-4(e^{15}+e^{26}), &\Omega^1_6=-4(e^{16}-e^{25})&\\[8pt]
\Omega^1_7=e^{13}-e^{17}+e^{24}-e^{28},& \Omega^1_8=e^{14}-e^{18}-e^{23}+e^{27},{}&\\[8pt]
\Omega^3_4=2(e^{12}-e^{34}),& \Omega^3_5=2(e^{35}+e^{46}),&{}\\[8pt]
\Omega^3_6=2(e^{36}-e^{45}), &\Omega^3_7=-(e^{37}+e^{48})&\\[8pt]
\Omega^3_8=2e^{34}-e^{38}-2e^{43}+e^{47},& \Omega^5_6=-8e^{56},&\\[8pt]
\Omega^5_7=2(e^{35}+e^{46}+e^{57}+e^{68}),& \Omega^5_8=2(-e^{36}+e^{45}+e^{58}-e^{67}),&{}\\[8pt]
\Omega^7_8=-2(e^{12}+e^{34}-e^{38}+e^{47}+e^{78})\,. &{}&
\end{array}
$$

\vskip.3cm

Therefore, Hol($\nabla^B$)=\SU(4).
}
\end{example}
\smallskip
\begin{example}
{\rm Let $\varphi^{j} = e^{2j-1} + i\, e^{2j}\,,\, j=1,\ldots ,4$ be complex $(1,0)$-forms satisfying
\begin{equation}\label{differentialcomplexforms}
d\varphi^1=0\,,\quad d\varphi^2=\varphi^{12}\,,\quad d\varphi^3=-\varphi^{13}\,,\quad
d\varphi^4=-\varphi^{23}\,.
\end{equation}
Hence
$$
\left\{
\begin{array}{ll}
de^1=0,\;\; de^2=0,\;\; de^3 =e^{13}-e^{24},\;\;& de^4 =e^{14}+e^{23},\;\; \\[6pt]
de^5 =-e^{15}+e^{26},\;\; de^6 =-e^{16}-e^{25},\;\;& de^7 =-e^{35}+e^{46},\;\; de^8 =-e^{36}-e^{45}\,.
\end{array}
\right.
$$
Let $\mathfrak{g}$ be the real $8$-dimensional Lie algebra whose
dual is spanned by $\{e^1,\ldots ,e^8\}$, and let $G$ be the
simply-connected Lie group whose Lie algebra is $\mathfrak{g}$. It
is immediate to see that $\mathfrak{g}$ is a $3$-step solvable but
not completely solvable Lie algebra. In fact, ${\mathfrak g}'
=\langle e_3,e_4,e_5,e_6,e_7,e_8\rangle\,, {\mathfrak g}''
=\langle e_7,e_8\rangle\,,{\mathfrak g}''' =\{0\}$. However, by
\cite{Nak} it turns out that $G$ has a uniform discrete subgroup
$\Gamma$, such that $M^8=\Gamma\backslash G$ is a compact
solvmanifold of dimension $8$.\newline Define an (integrable)
complex structure $J$ on $M^8$ by setting
$$
\begin{array}{l}
Je^1=-e^2,\;\;Je^2=e^1,\;\;Je^3=-e^4,\;\;Je^4=-e^3,\;\;\\[5pt]
Je^5=-e^6,\;\;Je^6=e^5,\;\;Je^7=-e^8,\;\;Je^8=e^7,
\end{array}
$$
and a $J$-Hermitian structure $g$ on $M$ by setting $g =\sum_{i=1}^8 e^i\otimes e^i$.\newline
Then $F=\sum_{j=1}^4e^{2j-1}\wedge e^{2j}$. Hence
$$
dF= 2\left(e^{134}-e^{156}\right)-e^{358}+e^{468}+e^{367}+e^{457}
$$
and consequently
$$
T= 2\left(e^{256}-e^{234}\right)-e^{467}+e^{357}+e^{458}+e^{368}\,.
$$
Therefore we obtain that the non zero connection and curvature forms of the Bismut connection are given respectively by
\begin{equation}\label{connectionexample4dim}
\begin{array}{llllll}
\omega^1_3=-e^3\,,&\omega^1_4=-e^4\,,&\omega^1_5=e^5\,,&\omega^1_6=e^6\,,&\omega^2_3=e^4\,,&\omega^2_4=-e^3\\[8pt]
\omega^2_5=-e^6\,,&\omega^2_6=e^5\,,&\omega^3_4=2e^2\,,&\omega^3_7=e^5\,,&\omega^3_8=e^6\,,&\omega^4_7=-e^6\\[8pt]
\omega^4_8=e^5\,,&\omega^5_6=-2e^2\,,&\omega^5_7=-e^3\,,&\omega^5_8=-e^4\,,&\omega^6_7=e^4\,,&\omega^6_8=-e^3
\end{array}
\end{equation}
and by
\begin{equation}\label{curvatureexample4dim}
\begin{array}{llll}
\Omega^1_2=2(e^{34}+e^{56}),&\Omega^1_3=-e^{24}-e^{13},&\Omega^1_4=e^{23}-e^{14},\!\!\!&\!\!\Omega^1_5=-e^{26}-e^{15}\\[8pt]
\Omega^1_6=e^{25}-e^{16},&\Omega^2_3=-e^{23}+e^{14},&\Omega^2_4=-e^{13}-e^{24},\!\!\!&\!\!\Omega^2_5=e^{16}-e^{25}\\[8pt]
\Omega^2_6=-e^{15}-e^{26},&\Omega^3_4=2(-e^{34}+e^{56}),&\Omega^3_7=-e^{26}-e^{15},\!\!\!&\!\!\Omega^3_8=e^{25}-e^{16}\\[8pt]
\Omega^4_7=-e^{25}+e^{16},&\Omega^4_8=-e^{26}-e^{15},&\Omega^5_6=2(-e^{56}+e^{34}),\!\!\!&\!\!\Omega^5_7=-e^{24}-e^{13}\\[8pt]
\Omega^5_8=e^{23}-e^{14},&\Omega^6_7=-e^{23}+e^{14},&\Omega^6_8=-e^{24}-e^{13},{}\!\!\!&\!\!\\[8pt]
\Omega^7_8=-2(e^{56}+e^{34})\,.{}&{}&{}&
\end{array}
\end{equation}
Hence the six $2$-forms
$$
\begin{array}{lll}
\Omega^1_2=2(e^{34}+e^{56})\,,\,&\Omega^1_3=-e^{24}-e^{13}\,,\,&\Omega^1_4=e^{23}-e^{14}\,,\\[8pt]
\Omega^1_5=-e^{26}-e^{15}\,,\,&\Omega^1_6=e^{25}-e^{16}\,,\,&\Omega^3_4=2(-e^{34}+e^{56})
\end{array}
$$
are linearly independent. Therefore, $\dim (\mathfrak{hol}$ $(\nabla^B))\geq 6$. \newline
We need to compute the covariant derivative of the curvature.
By the above expression for the connection forms, we get:
\begin{eqnarray*}
\nabla^B e^1 &=& e^3\otimes e^3 + e^4\otimes e^4 - e^5\otimes e^5 - e^6\otimes e^6\\
\nabla^B e^2 &=& -e^4\otimes e^3 + e^3\otimes e^4 + e^6\otimes e^5 - e^5\otimes e^6\\
\nabla^B e^3 &=& -e^3\otimes e^1 + e^4\otimes e^2 -2\, e^2\otimes e^4 - e^5\otimes e^7 - e^6\otimes e^8\\
\nabla^B e^4 &=& -e^4\otimes e^1 - e^3\otimes e^2 +2\, e^2\otimes e^3 + e^6\otimes e^7 - e^5\otimes e^8\\
\nabla^B e^5 &=& e^5\otimes e^1 - e^6\otimes e^2 +2\, e^2\otimes e^6 + e^3\otimes e^7 + e^4\otimes e^8\\
\nabla^B e^6 &=& e^6\otimes e^1 + e^5\otimes e^2 - 2\, e^2\otimes e^5 - e^4\otimes e^7 + e^3\otimes e^8\\
\nabla^B e^7 &=& e^5\otimes e^3 - e^6\otimes e^4 - e^3\otimes e^5 + e^4\otimes e^6\\
\nabla^B e^8 &=& e^6\otimes e^3 + e^5\otimes e^4 - e^4\otimes e^5 - e^3\otimes e^6\,.
\end{eqnarray*}
By the above formulas, a straightforward computation, taking into account \eqref{curvatureexample4dim},  gives
\begin{eqnarray*}
\nabla^B_{e_3}\Omega^1_2 &=& 2\left(-e^{14}+e^{23}-e^{67}+e^{58}\right)\\
\nabla^B_{e_4}\Omega^1_2 &=& 2\left(e^{13}+e^{24}-e^{68}-e^{57}\right)\\
\nabla^B_{e_5}\Omega^1_2 &=& 2\left(e^{47}-e^{38}+e^{16}-e^{25}\right)\\
\nabla^B_{e_6}\Omega^1_2 &=& 2\left(e^{48}+e^{37}-e^{26}-e^{15}\right)\\
\nabla^B_{e_6}\Omega^1_4 &=& -e^{35}-e^{28}-e^{46}-e^{17}\\
\nabla^B_{e_4}\Omega^1_5 &=& -e^{18}+e^{27}+e^{36}-e^{45}\\
\nabla^B_{e_3}\Omega^1_6 &=& e^{45}+e^{27}-e^{36}-e^{18}\\
\nabla^B_{e_4}\Omega^1_6 &=& -e^{35}+e^{28}-e^{46}+e^{17}\\
\nabla^B_{e_5}\Omega^1_6 &=& 2(e^{56}-e^{12})\,.
\end{eqnarray*}
Therefore, we have proved the following
\begin{proposition} The holomorphic parallelizable solvmanifold $M=\Gamma\backslash G$, where
$\{\varphi^1,\ldots , \varphi^4\}$ are the complex $(1,0)$-forms satisfying
\eqref{differentialcomplexforms}, it is endowed with an Hermitian metric with
K\"ahler form
$$
F=\frac{i}{2}\sum_{j=1}^4\varphi^j\wedge \overline{\varphi}^j
$$
satisfying $dF^3=0$ and the infinitesimal holonomy of the Bismut connection is ${\mathfrak su}(4)$.
\end{proposition}
}
\end{example}

\medskip

\noindent {\bf Acknowledgments.} This work has been partially
supported through grant MCyT (Spain) MTM2005-08757-C04-02, by the Project MIUR
``Geometric Properties of Real and Complex Manifolds'' and by GNSAGA of INdAM.

\smallskip

{\small

\end{document}